\documentclass[times,review,sortcompress]{elsarticle}
\usepackage[active]{srcltx}
\usepackage{hyperref}
\usepackage{geometry}
\usepackage{bbm}
\usepackage{mathrsfs}
\usepackage{tipa}
\usepackage{amsthm}
\usepackage{amssymb}
\usepackage{stackrel}
\usepackage{amsmath,amscd}
\usepackage[mathscr]{eucal}
\usepackage{mathrsfs}
\usepackage[active]{srcltx}
\usepackage{amsfonts,amssymb}
\usepackage{amsthm,amsmath}
\usepackage[all]{xy}
\usepackage{geometry}
\usepackage{rotating}
\usepackage{lscape}
\usepackage{longtable}
\usepackage{booktabs}
\usepackage[active]{srcltx}
\usepackage[centerlast]{caption}
\usepackage{wallpaper}
\usepackage{graphics}
\usepackage{graphicx}
\usepackage{floatrow}
\usepackage{caption}
\usepackage{boxedminipage}
\usepackage{shadow}
\newtheorem{dl}{Theorem}[section]
\newtheorem{tl}[dl]{Corollary}
\newtheorem{yl}[dl]{Lemma}
\newtheorem{dy}[dl]{Definition}
\newtheorem{lz}[dl]{Example}
\newtheorem{wenti}[dl]{Research problem}

\newtheorem{remark}[dl]{Remark}

\numberwithin{equation}{section}
\newproof{pf}{Proof}
\newproof{pot333}{Proof of Identity \eqref{generaltrans-7}}
\newproof{pot444}{Proof of Identity \eqref{generaltrans-7-7-0}}
\newproof{pot555}{Proof of Identity \eqref{4.1}}
\newproof{pot666}{Proof of Theorem \ref{mainthm-good}}
\newproof{pot111}{Proof of Theorem \ref{mainthm-qidentities}}
\newproof{pot222}{Proof of Theorem \ref{baileylemma}}
\newcommand{\poq}[2]{(#1;q)_{#2}}

\def\qed{\hfill \rule{4pt}{7pt}}
\def\pf{\noindent {\it Proof.} }
\newcommand{\poqq}[2]{(#1;q^2)_{#2}}

\begin{document}
\title{Further results on Andrews--Yee's two identities for mock theta functions $\omega(z;q)$ and $v(z;q)$}
\author{Jin Wang\fnref{a,fn1}}
\fntext[fn1]{E-mail address: jinwang2016@yahoo.com}
\address[a]{Department of Mathematics, Soochow University, SuZhou 215006, P.R.China}
\author{Xinrong Ma\fnref{b,fn2,fn3}}
\fntext[fn2]{This  work was supported by NSFC grant No. 11471237}
\fntext[fn3]{Corresponding author. E-mail address: xrma@suda.edu.cn.}
\address[b]{Department of Mathematics, Soochow University, SuZhou 215006, P.R.China}
\markboth{J. Wang and  X. Ma}{On Andrews--Yee's   identities of  mock theta functions}
\begin{abstract}
In this paper, by the method of comparing coefficients and the inverse technique, we establish the corresponding
variate forms of two  identities of  Andrews and Yee for mock theta functions, as well as a few allied but unusual
$q$-series identities. Among includes a new Bailey pair from which a product formula of two ${}_2\phi_1$ series is
derived.  Further, we focus on two finite $q$-series summations arising from  Andrews and Yee's mock theta function
identities  and expound some recurrence relations  and transformation formulas behind them.
\end{abstract}
\maketitle
\vspace{20pt}
\parskip 7pt
\baselineskip 16pt
\section{Introduction}
In the very recent  paper  \cite{andrews}, G.E. Andrews and A.J. Yee   showed the following  $q$-series identities.
\begin{yl}[\mbox{\cite[Theorem 1, Eq.(6)/Eq.(7); Eq.(8)]{andrews}}]\label{andrews-yee}For any complex numbers $q:|q|<1$
and variable $z$, the following identities are valid.
\begin{align}
\displaystyle\sum_{n=0}^\infty
\frac{z^nq^{2n^2+2n+1}}{\poqq{q,zq}{n+1}}&=\sum_{n=1}^\infty\frac{z^{n-1}q^{n}}{\poqq{q}{n}},\label{4.1}\\
\sum_{n=0}^\infty\,q^n\poq{-zq^{n+1}}{n}
\poqq{zq^{2n+2}}{\infty}&=\sum_{n=0}^{\infty}\frac{z^{n}q^{n^2+n}}{\poqq{q}{n+1}},\label{generaltrans-7-7}\\
\sum_{n=1}^\infty
\frac{q^n}{\poq{zq^n}{n+1}\poqq{zq^{2n+2}}{\infty}}&=\sum_{n=1}^\infty\frac{z^{n-1}q^{n}}{\poqq{q}{n}}.\label{generaltrans-7-77}
\end{align}
\end{yl}
According to G.E. Andrews and A.J. Yee,  these $q$-series identities are closely related with
two mock theta
functions $\omega(z;q)$ and $v(z;q)$ which were first introduced by G.E. Andrews \cite{andrews0} as follows:
\begin{align*}
\omega(z;q):=\sum_{n=0}^\infty \frac{z^nq^{2n^2+2n}}{\poqq{q}{n+1}\poqq{zq}{n+1}}\,\, \mbox{and}\,\,
v(z;q):=\sum_{n=0}^{\infty}\frac{q^{n^2+n}}{\poqq{-zq}{n+1}}.
\end{align*}
As is well known to us,  $\omega(z;q)$ and $v(z;q)$ go back to
 the third order mock theta functions $\omega(q)=\omega(1;q)$ and $v(q)=\omega(1;q)$ due to G.N. Watson \cite{watson}
and afterward rediscovered in Ramanujan's lost notebook \cite{ramanujan}. Perhaps the most important thing is that
 in their paper \cite{dixit}, G.E. Andrews, A. Dixit and A.J. Yee discovered that  $\omega(q)$ and $v(q)$ serve as
 enumerative functions for special integer partitions. Regarding this together with more recent progress, we refer the
 reader to
 \cite{dixit} and \cite{newadded} of Y.S. Choi.

  In the present paper, we shall generalize the above  Andrews--Yee identities  \eqref{generaltrans-7-7} and
  \eqref{generaltrans-7-77} to the following variate (i.e. $y$) forms and then prove \eqref{4.1} in a different way.
\begin{dl}\label{mainthm} For arbitrary complex numbers $q, y$  such that $|q|, |y|<1$, the following identities are
valid.
\begin{align}
\sum_{n=0}^\infty\,y^n\poq{-zq^{n+1}}{n}\poqq{zq^{2n+2}}{\infty}&=\sum_{n=0}^{\infty}z^nq^{n^2+n}\frac{\poq{-y}{n}}{\poq{yq^n}{n+1}}
\sum_{k=0}^{n} q^{2k}\frac{(y/q;q){}_{2k}}
{ \left(q^2,y^2;q^2\right){}_{k}},\label{generaltrans-7-7-0}\\
\sum_{n=1}^\infty
\frac{y^{n-1}}{\poq{zq^n}{n+1}\poqq{zq^{2n+2}}{\infty}}&=\sum_{n=0}^\infty\,(qz)^n\frac{\poq{-y}{n}}{\poq{yq^n}{n+1}}\sum
_{k=0}^{n} q^k\frac{\poq{y/q}{2k}}{\poqq{q^2,y^2}{k}}.\label{generaltrans-7}
\end{align}
\end{dl}
In the process of  proving Theorem \ref{mainthm}, we obtain some special $q$-series identities. To our best knowledge,
they are new and of interest  deserving separate theorems below.
\begin{dl}\label{mainthm-qidentities}  For integer $n\geq 0$ and complex numbers  $y$ being not of the form $-q^{-n}$,
the following identities are valid.
\begin{align}
\sum_{k=0}^n q^k\frac{\poq{yq^{n}}{k}}{\poqq{q^2}{k}}&=\poq{-y}{n}\sum _{{\it
k}=0}^{n}q^k\frac{\poq{y/q}{2k}}{\poqq{q^2,y^2}{k}} ,\label{deduce111}\\
\sum_{k=0}^n\bigg[\genfrac{}{}{0pt}{}{n}{k}\bigg]_qy^kq^{k^2/2-k/2}\frac{\poq{q}{k}}{\poq{yq^{n}}{k+1}}&=\frac{(q^2,y^2;q^2){}_n
}{(y;q){}_{2n+1}}
\sum_{k=0}^{n} q^{2k} \frac{(y/q;q){}_{2k}}
{ \left(q^2,y^2;q^2\right){}_{k}}\label{deduce222}.
\end{align}
In particular, we have
\begin{align}
\sum_{k=0}^\infty
y^kq^{k^2/2-k/2}=(q^2;q^2){}_\infty(-y;q){}_{\infty}\,{}_{2}\phi_1\bigg[\genfrac{}{}{0pt}{}{y/q,y}{y^2};q^2,q^2\bigg]\label{thetanew}
\end{align}
while, for any integer $r\geq 0$,
\begin{align}
\poq{-q}{2r}\sum_{k=0}^\infty q^{k^2/2+k/2+2rk}=\frac{(q^2;q^2){}_\infty}{\poqq{q}{\infty}}\sum_{k=0}^{\infty} q^{2k}
\frac{(q^{2r};q){}_{2k}}
{ \left(q^2,q^{4r+2};q^2\right){}_{k}}.\label{thetanew-added}
\end{align}
\end{dl}

Especially noteworthy is that  \eqref{deduce222} offers a new Bailey pair and  \eqref{thetanew-added} is also different from one given by G.E. Andrews and S.O. Warnaar in
\cite{andrewswarnaar2007} via the Bailey transformation. Actually, as a combination of \eqref{deduce222} and the classical Bailey lemma, we can set up

\begin{dl}\label{baileylemma} With the same conditions as Theorem \ref{mainthm-qidentities}. Then for any  $|yq/ab|<1$,
we have
\begin{align}
\sum_{n=0}^{\infty}\frac{\poq{a,b}{n}}
{\poq{yq/a,yq/b}{n}}\left(\frac{y^2}{ab}\right)^nq^{n^2/2+n/2}\label{slater-result}\\
= \frac{\poq{yq,yq/ab}{\infty}} {\poq{yq/a,yq/b}{\infty}}
\sum_{n=0}^{\infty}\frac{(-q,-y,a,b;q){}_n }{(yq;q){}_{2n}}\left(\frac{yq}{ab}\right)^n
\sum_{k=0}^{n} \frac{q^{2k}(y/q;q){}_{2k}}
{\left(q^2,y^2;q^2\right){}_{k}}.\nonumber
 \end{align}
 \end{dl}

It comes as a surprise that  a product formula  of two ${}_2\phi_1$ series emerges from  \eqref{deduce222}.  It is
similar to but different from the
existing results such as various $q$-analogues of the Clausen formula \cite[Section 8.8]{10}.
\begin{dl}\label{mainthm-good} For any complex numbers $a$ and $b$ being not of the form $-q^{-m}$, integer $m\geq 0$,
it holds
\begin{align}
{}_{2}\phi_1\bigg[\genfrac{}{}{0pt}{}{a,a/q}{a^2};q^2,q^2\bigg]{}_{2}\phi_1\bigg[\genfrac{}{}{0pt}{}{b,b/q}{b^2};q^2,q^2\bigg]
=\frac{\poqq{q}{\infty}}{\poqq{q^2}{\infty}}{}_{4}\phi_3\bigg[\genfrac{}{}{0pt}{}{(ab)^{1/2},-(ab)^{1/2},(ab/q)^{1/2},-(ab/q)^{1/2}}{-a,-b,ab/q};q,q\bigg].\label{thetanewidentity}
\end{align}
\end{dl}
 Before proceeding, let us  give several remarks on notation.  Throughout this paper, we adopt the standard notation and
 terminology for  $q$-series from the book  \cite{10}.  As customary, the $q$-shifted factorials of complex variable $x$
 with the base $q$ are  given by
\begin{eqnarray*}
(x;q)_\infty
:=\prod_{n=0}^{\infty}(1-xq^n)\qquad\mbox{and}\quad (x;q)_n:=\frac{(x;q)_\infty}{(xq^n;q)_\infty}
\end{eqnarray*}
for all integers $n$. For integers $m\geq 1$, we  employ
 the multi-parameter compact notation
\[(a_1,a_2,\ldots,a_m;q)_n:=(a_1;q)_\infty (a_2;q)_n\ldots (a_m;q)_n.\]
Also, the $q$-binomial coefficients $\bigg[\genfrac{}{}{0pt}{}{n}{k}\bigg]_q$ and the ${}_{r+1}\phi_r$ series with the
base $q$ and the argument $x$ are respectively defined
\begin{align*}
\bigg[\genfrac{}{}{0pt}{}{n}{k}\bigg]_q&:=\frac{\poq {q}{n}}{\poq
{q}{n-k}\poq{q}{k}},\\
{}_{r+1}\phi_r\bigg[\genfrac{}{}{0pt}{}{a_1,a_2,\cdots,a_{r+1}}{b_1,b_2,\cdots,b_{r}};q,x\bigg]&:=\sum _{n=0} ^{\infty}
\frac{\poq{a_1,a_2,\ldots,a_{r+1}}{n}}{\poq{q,b_1,b_2,\ldots,b_r}{n}}x^{n}.
\end{align*}
To easy notation, we write $\tau(n)$ for $(-1)^nq^{n(n-1)/2}$.
In addition, given any power series $f(x)=\sum_{n\geq 0}a_nx^n$ in variable $x$,  we define the coefficient functional
$$\boldsymbol\lbrack x^n\boldsymbol\rbrack f(x):=a_n\,\,\mbox{and}\,\, a_0=f(0).$$

Our paper is organized as follows. In Section 2, Theorems \ref{mainthm}--\ref{mainthm-good}, as well as Identity
\eqref{4.1}  will be proved.  In
Section 3, we shall study two finite $q$-series sums related to these mock theta function identities. As applications,
 some interesting $q$-series identities will be presented in these two sections. Among includes a new Bailey pair from
 which a product formula of two ${}_2\phi_1$ series is derived.  In Section 4, we conclude our
paper with some remarks on $\omega(z;q)$ and $v(z;q)$.

\section{Proofs of main theorems}
Our argument rests on the following  transformations of terminating  ${}_3\phi_2$ series.
\begin{yl}[\mbox{\cite[(III.11)/(III.13)]{10}}]\label{yl2.1}
\begin{align}
&{_3\phi_2}\bigg[\genfrac{}{}{0pt}{}{q^{-n},b,c}{d,e};q,q\bigg]=\frac{\poq{de/bc}{n}}{\poq{e}{n}}\bigg(\frac{bc}{d}\bigg)^n
{_3\phi_2}\bigg[\genfrac{}{}{0pt}{}{q^{-n},d/b,d/c}{d,de/bc};q,q\bigg],\label{generaltrans-new-1}\\
&{_3\phi_2}\bigg[\genfrac{}{}{0pt}{}{q^{-n},b,c}{d,e};q,\frac{deq^n}{bc}\bigg]=\frac{\poq{e/c}{n}}{\poq{e}{n}}
{_3\phi_2}\bigg[\genfrac{}{}{0pt}{}{q^{-n},c,d/b}{d,cq^{1-n}/e};q,q\bigg].\label{generaltrans-new}
\end{align}
\end{yl}
\subsection{Proof of  Theorem \ref{mainthm}}
Next is the complete proof of  Theorem \ref{mainthm}, which is composed of  the proofs of  \eqref{generaltrans-7-7-0}
and \eqref{generaltrans-7}.
\begin{pot444}  To establish \eqref{generaltrans-7-7-0}, we temporarily
  assume that
\begin{align}
\sum_{n=0}^\infty\,y^n\poq{-zq^{n+1}}{n}\poqq{zq^{2n+2}}{\infty}=\sum_{n=0}^{\infty}f_n(y)z^n.\label{generaltrans-7-7-00}
\end{align}
It is easy to see that
$$
\poq{-zq^{n+1}}{n}\poqq{-zq^{2n+2}}{\infty}=\frac{\poq{-zq^{n+1}}{\infty}}{\poqq{-zq^{2n+1}}{\infty}}.
$$
By the $q$-binomial theorem \cite[(II.3)]{10}, we have
\begin{align*}
\frac{\poq{-zq^{n+1}}{\infty}}{\poqq{-zq^{2n+1}}{\infty}}&=\sum_{i=0}^\infty\frac{\tau(i)}{\poq{q}{i}}
(-zq^{n+1})^i\sum_{j=0}^\infty\frac{(-zq^{2n+1})^j}{\poqq{q^2}{j}}.
\end{align*}
 Therefore, by equating the coefficients of $z^m$ on both sides, to show \eqref{generaltrans-7-7-0} is to find
\begin{align*}
f_m(y)&=\sum_{n=0}^\infty
y^n\sum_{i+j=m}(-1)^i\frac{\tau(i)q^{ni+i}}{\poq{q}{i}}\frac{(-1)^j(q^{2n+1})^j}{\poqq{q^2}{j}}\\
&=(-1)^mq^m\sum_{i+j=m}\frac{\tau(i)}{\poq{q}{i}\poqq{q^2}{j}}\sum_{n=0}^\infty q^{n(i+2j)}y^n.
\\
&=(-1)^mq^m\sum_{i+j=m}\frac{\tau(i)}{\poq{q}{i}\poqq{q^2}{j}}\frac{1}{1-yq^{m+j}}.
\end{align*}
Next, we multiply both sides by $(-1)^mq^{-m}(1-yq^m)\poq{q}{m}$ and reformulate it as
\begin{align}
(-1)^mq^{-m}(1-yq^m)\poq{q}{m}f_m(y)=\sum_{j=0}^m\bigg[\genfrac{}{}{0pt}{}{m}{m-j}\bigg]_q\tau(m-j)\frac{\poq{yq^{m}}{j}}{\poq{-q,yq^{1+m}}{j}}.\label{generaltrans-8888-ii-0}
\end{align}
In this form, we need to examine the term
\begin{align*}
\bigg[\genfrac{}{}{0pt}{}{m}{m-j}\bigg]_q\tau(m-j)&=\frac{\poq{q^{-m}}{m-j}}{\poq{q}{m-j}}q^{m(m-j)}
=\frac{\poq{q^{-m}}{m}}{\poq{q}{m}}\frac{\poq{q^{-m}}{j}}{\poq{q}{j}}q^{m(m-j)+(m+1)j}\\
&=\tau(m)q^j\frac{\poq{q^{-m}}{j}}{\poq{q}{j}},
\end{align*}
wherein we have used the basic relations (cf.\cite[(I.42)]{10})
\begin{equation}
\bigg[\genfrac{}{}{0pt}{}{n}{k}\bigg]_q\tau(k)=\frac{\poq{q^{-n}}{k}}{\poq{q}{k}}q^{nk}\label{basic1}
\end{equation}
and (cf.\cite[(I.11)]{10})
\begin{equation}\frac{\poq{a}{n-k}}{\poq{b}{n-k}}=\frac{\poq{a}{n}}{\poq{b}{n}}
\frac{\poq{q^{1-n}/b}{k}}{\poq{q^{1-n}/a}{k}}(b/a)^k.\label{basic2}\end{equation}
Accordingly, \eqref{generaltrans-8888-ii-0} reduces to
\begin{equation}
(1-yq^m)(-1)^mq^{-m}\poq{q}{m}f_m(y)=\tau(m){_3\phi_2}\bigg[\genfrac{}{}{0pt}{}{q^{-m},yq^{m},0}{-q,yq^{1+m}};q,q\bigg].\label{generaltrans-8888-iii-00}
\end{equation}
It is now clear that we can employ the transformation \eqref{generaltrans-new-1} under the simultaneous parametric
specializations
$$(n,b,c,d,e)=(m,yq^{m},0,yq^{1+m},-q)$$ to the ${}_3\phi_2$ series on the right-hand of
\eqref{generaltrans-8888-iii-00},
arriving at
\begin{equation*}
{_3\phi_2}\bigg[\genfrac{}{}{0pt}{}{q^{-m},yq^{m},0}{-q,yq^{1+m}};q,q\bigg]=\frac{q^{m^2/2+m/2}}{\poq{-q}{m}}{_2\phi_1}\bigg[\genfrac{}{}{0pt}{}{q^{-m},q}{yq^{m+1}};q,-yq^{m}\bigg].
\end{equation*}
Making use of this fact, we see that  \eqref{generaltrans-8888-iii-00} is equivalent to
\begin{equation*}(-1)^mq^{-m}(1-yq^m)\poq{q}{m}f_m(y)=\frac{(-1)^mq^{m^2}}{\poq{-q}{m}}
\sum_{k=0}^m\bigg[\genfrac{}{}{0pt}{}{m}{k}\bigg]_qy^kq^{k^2/2-k/2}\frac{\poq{q}{k}}{\poq{yq^{m+1}}{k}}
,
\end{equation*}
or equivalently,
\begin{equation}q^{-m^2-m}\poqq{q^2}{m}f_m(y)=\sum_{k=0}^m\bigg[\genfrac{}{}{0pt}{}{m}{k}\bigg]_qy^kq^{k^2/2-k/2}\frac{\poq{q}{k}}{\poq{yq^{m}}{k+1}}.
\label{generaltrans-888888-999}
\end{equation}
Now we can solve  \eqref{generaltrans-888888-999} for $f_m(y)$. To that end,  we  define
\begin{align*}
S_m(k):=\bigg[\genfrac{}{}{0pt}{}{m}{k}\bigg]_qy^kq^{k^2/2-k/2}\frac{\poq{q}{k}}{\poq{yq^{m}}{k+1}}\quad\mbox{and}\quad
T_m:=\sum_{k=0}^mS_m(k).
\end{align*}
Next, by the WZ method (see \cite{aeqb} for further information), it holds
\begin{align*}S_{m+1}(k)-\frac{(1-q^{2m+2}) (1-y^2 q^{2m})
   }{(1-y q^{2 m+1}) (1-y q^{2
   m+2})}S_{m}(k)=G(m,k)-G(m,k-1),
\end{align*}
where
\[G(m,k):=\frac{
  y q^{k+m+1}+ q^{k+1}-
   q^{m+1}-1}{\left(1-y q^{2 m+1}\right) \left(1-y q^{2 m+2}\right)
  }\times\frac{\poq{q}{k+1}}{\left(yq^{m+1};q\right){}_{k+1}}
  \bigg[\genfrac{}{}{0pt}{}{m+1}{k+1}\bigg]_qy^{k+1} q^{\frac{k^2}{2}-\frac{k}{2}+m}
.\]
By telescopying, we find that $T_0=1/(1-y)$ and
\begin{align}T_{m}-\frac{(1-q^{2m}) (1-y^2 q^{2m-2})
   }{(1-y q^{2 m-1}) (1-y q^{2m})}T_{m-1}=\frac{q^{2 m-1} (q-y)}{\left(1-y q^{2 m-1}\right) \left(1-y q^{2
   m}\right)}.\label{t_m}
\end{align}
Solving \eqref{t_m}, we obtain
\begin{align*}
T_m=\frac{(q^2,y^2;q^2){}_m }{(y;q){}_{2m+1}}
\sum_{k=0}^{m} \frac{q^{2k}(y/q;q){}_{2k}}
{ \left(q^2,y^2;q^2\right){}_{k}},
\end{align*}
 which in  turn yields
\begin{align}f_m(y)=q^{m^2+m}\frac{(y^2;q^2){}_m }{(y;q){}_{2m+1}}
\sum_{k=0}^{m} \frac{q^{2k}(y/q;q){}_{2k}}
{ \left(q^2,y^2;q^2\right){}_{k}}.\label{results-1}
\end{align}
A direct substitution of \eqref{results-1} back to  \eqref{generaltrans-7-7-00} completes the proof of
\eqref{generaltrans-7-7-0}.
\qed
\end{pot444}

 Now we proceed to the proof of Identity  \eqref{generaltrans-7}.
\begin{pot333} As previously, in order to show \eqref{generaltrans-7},  we assume for the moment that
\begin{align}
\sum_{n=1}^\infty
\frac{y^n}{\poq{zq^n}{n+1}\poqq{zq^{2n+2}}{\infty}}=\sum_{n=0}^\infty\,g_n(y)z^n.\label{generaltrans-70}
\end{align}
Observe that
\begin{align*}
\frac{1}{\poq{zq^n}{n+1}\poqq{zq^{2n+2}}{\infty}}=\frac{\poqq{zq^{2n+1}}{\infty}}{\poq{zq^n}{\infty}}.
\end{align*}
So, by comparing coefficients, we see that for any integer $m\geq 0$,
\begin{align*}
 g_m(y)=\sum_{n=1}^\infty\, y^n[z^m]\frac{\poqq{zq^{2n+1}}{\infty}}{\poq{zq^n}{\infty}}.
\end{align*} It remains to find $g_m(y)$. For this, we invoke the $q$-binomial theorem \cite[(II.3)]{10} to expand
\begin{align*}
\frac{\poqq{zq^{2n+1}}{\infty}}{\poq{zq^n}{\infty}}=\sum_{i=0}^\infty\frac{(zq^n)^i}{\poq{q}{i}}
\sum_{j=0}^\infty\frac{(-1)^jq^{j^2-j}(zq^{2n+1})^j}{\poqq{q^2}{j}}.
\end{align*}
With this given, we need to compute
\begin{align}
g_m(y)&=\sum_{n=1}^\infty
y^n\sum_{i+j=m}\frac{q^{ni}}{\poq{q}{i}}\frac{(-1)^jq^{j^2-j}(q^{2n+1})^j}{\poqq{q^2}{j}}\nonumber\\
&=\sum_{i+j=m}\frac{(-1)^jq^{j^2}}{\poq{q}{i}\poqq{q^2}{j}}\sum_{n=1}^\infty
(y q^{i+2j})^n\nonumber\\
&=\frac{yq^m}{1-yq^m}\sum_{i+j=m}\frac{(-1)^jq^{j^2+j}}{\poq{q}{i}\poq{q}{j}\poq{-q}{j}}\frac{\poq{yq^m}{j}}{\poq{yq^{m+1}}{j}}.\label{bbb}
\end{align}
Recasting the rightmost sum in \eqref{bbb}  in terms of the ${}_{3}\phi_2$ series, we see that
\begin{equation}
(1-yq^{m})\poq{q}{m}g_m(y)=yq^m\lim_{b\to\infty}{_3\phi_2}\bigg[\genfrac{}{}{0pt}{}{q^{-m},yq^{m},b}{-q,yq^{m+1}};q,-\frac{q^{2+m}}{b}\bigg].\label{generaltrans-88888-x}
\end{equation}
At this stage, we apply the  transformation  \eqref{generaltrans-new}   under the simultaneous parametric specialization
$$(n,b,c,d,e)=(m,b,yq^{m},-q,yq^{m+1})$$
to the ${}_3\phi_2$ series on the right-hand side of  \eqref{generaltrans-88888-x}, thereby obtaining
\begin{equation*}
{_3\phi_2}\bigg[\genfrac{}{}{0pt}{}{q^{-m},yq^{m},b}{-q,yq^{1+m}};q,-\frac{q^{2+m}}{b}\bigg]
=\frac{\poq{q}{m}}{\poq{yq^{1+m}}{m}}{_3\phi_2}\bigg[\genfrac{}{}{0pt}{}{q^{-m},yq^{m},-q/b}{-q,q^{-m}};q,q\bigg].
\end{equation*}
With the aid of this relation, it is easy to see that \eqref{generaltrans-88888-x} is equivalent to
\begin{align*}(1-yq^{m})\poq{q}{m}g_m(y)=yq^m\frac{\poq{q}{m}}{\poq{yq^{1+m}}{m}}\sum_{k=0}^m\frac{\poq{yq^{m}}{k}}{\poqq{q^2}{k}}q^k.
\end{align*}
After simplification, it becomes
\begin{align}
\frac{\poq{yq^{m}}{m+1}}{yq^m}g_m(y)=\sum_{k=0}^m\frac{\poq{yq^{m}}{k}}{\poqq{q^2}{k}}q^k.\label{definition}
\end{align}
In order to find  $g_m(y)$, we invoke  the WZ method  and define
\begin{align*}
S_m(k):=\frac{\poq{yq^{m}}{k}}{\poqq{q^2}{k}}q^k \quad\mbox{and}\quad T_m:=\sum_{k=0}^mS_m(k).
\end{align*}
In both notation, we suppress the dependence on $y$. Next, by the WZ method, we find
\begin{align*}S_{m+1}(k)-(1+yq^{m})S_{m}(k)=G(m,k)-G(m,k-1),
\end{align*}
where
\[G(m,k):=-y q^{m}\frac{\poq{yq^{m+1}}{k}}{\poqq{q^2}{k}}.\]
By telescopying, it is without difficult to check that $T_m$  satisfies
\begin{align*}T_{m+1}-(1+yq^{m})T_{m}
=\frac{q^{m} (q-y) \poq{yq^{m+1}}{m}}{\poqq{q^2}{m+1}}.
\end{align*}
Solving  this recurrence relation,
we get
\begin{align*}
T_m=\poq{-y}{m}\sum _{k=0}^{m} q^k\frac{\poq{y/q}{2k}}{\poqq{q^2,y^2}{k}}.
\end{align*}
 And then substituting this for $T_m$ of \eqref{definition} and solving the resulted for $g_m(y)$, we immediately obtain
\begin{align}
g_m(y)=yq^m\frac{\poq{-y}{m}}{\poq{yq^m}{m+1}}\sum _{k=0}^{m} q^k\frac{\poq{y/q}{2k}}{\poqq{q^2,y^2}{k}}.\label{fresult}
\end{align}
A direct substitution of \eqref{fresult} back to  \eqref{generaltrans-70}  gives the complete proof of the identity.
\qed
\end{pot333}

\subsection{Proof of Theorem \ref{mainthm-qidentities}}
Now we are ready to present a short proof of Theorem  \ref{mainthm-qidentities}.
\begin{pot111} It is clear  from the above argument that \eqref{deduce111} is obtained by combining  \eqref{definition}
with \eqref{fresult} while  \eqref{deduce222} is the consequence of substituting  \eqref{results-1} into
\eqref{generaltrans-888888-999}. On taking the limitation of \eqref{deduce222} as  $n$ tends to infinite, we arrive at
\eqref{thetanew}. It in turn yields  \eqref{thetanew-added} under $y=q^{2r+1}$. \qed
\end{pot111}

\subsection{Proofs of  Theorems \ref{baileylemma} and \ref{mainthm-good}}
Now we proceed to show  Theorem \ref{baileylemma} via the use of \eqref{deduce222}. To this purpose, we need to recall
the concept of Bailey pair relative to $t$ which can be found in the book \cite{andrews2} by G.E. Andrews.
\begin{dy}\label{baileydef}
A Bailey pair relative to $t$ is conveniently  defined to be a pair of sequences $\{\alpha_n(t)\}_{n\geq 0}$ and
$\{\beta_n(t)\}_{n\geq 0},$ denoted by  $(\alpha_n(t),\beta_n(t))$,  satisfying
\begin{eqnarray}
\beta_n(t)=\sum^{n}_{k=0} \frac{\alpha_k(t)}
{\poq{q}{n-k}\poq{tq}{n+k}}.\label{defnew}
 \end{eqnarray}
\end{dy}
Closely related to Bailey pairs is  the well-known  Bailey lemma, which first appeared in \cite[Eq.(3.1)]{bailey} of
W.N. Bailey. See also \cite[Eq.(3.4.9)]{sla} or \cite[Eq.(3.33)]{andrews2} for reference.
\begin{yl}[Bailey  lemma]
\begin{align}
\sum_{n=0}^{\infty}\frac{\poq{a,b}{n}}
{\poq{tq/a,tq/b}{n}}\left(\frac{tq}{ab}\right)^n\alpha_n(t)\label{slater}
= \frac{\poq{tq,tq/ab}{\infty}} {\poq{tq/a,tq/b}{\infty}}
\sum_{n=0}^{\infty}\poq{a,b}{n}\left(\frac{tq}{ab}\right)^n\beta_n(t)
 \end{align}
provided that all the
relevant infinite series absolutely convergent.
  \end{yl}

For a good survey on Bailey pairs and various Bailey lemmas, the reader is referred to \cite{warnaar0} due to S.O.
Warnaar. Especially noteworthy  here is that, as demonstrated by  G.E. Andrews and B.C. Berndt in \cite[\S 11.5]{ber}
and \cite[Chap. 5]{ber-1},   the  Bailey  lemma is a basic tool to the study of Ramanujan's mock theta function
identities.

\begin{pot222}
In view of Definition \ref{baileydef}, we now reformulate \eqref{deduce222} in the form
\begin{eqnarray*}
\beta_n(y)=\sum^{n}_{k=0}\frac{\alpha_k(y)}
{\poq{q}{n-k}\poq{yq}{n+k}},
 \end{eqnarray*}
 where $\alpha_n(y)$ and $\beta_n(y)$ are, respectively, given by
 \begin{align}\label{defnew-2}
  \begin{cases}
   \,\alpha_n(y)=\displaystyle  y^nq^{n^2/2-n/2}\\[2mm] \displaystyle
   \, \beta_n(y)=\frac{(-q,-y;q){}_n }{(yq;q){}_{2n}}
\sum_{k=0}^{n} \frac{q^{2k}(y/q;q){}_{2k}}
{\left(q^2,y^2;q^2\right){}_{k}}.
  \end{cases}
 \end{align}
Evidently, $(\alpha_n(y),\beta_n(y))$ is such a Bailey pair that  reduces \eqref{slater} to
\begin{align*}
\sum_{n=0}^{\infty}\frac{\poq{a,b}{n}}
{\poq{yq/a,yq/b}{n}}\left(\frac{y^2q}{ab}\right)^nq^{n^2/2-n/2}\\
= \frac{\poq{yq,yq/ab}{\infty}} {\poq{yq/a,yq/b}{\infty}}
\sum_{n=0}^{\infty}\poq{a,b}{n}\left(\frac{yq}{ab}\right)^n\frac{(-q,-y;q){}_n }{(yq;q){}_{2n}}
\sum_{k=0}^{n} \frac{q^{2k}(y/q;q){}_{2k}}
{\left(q^2,y^2;q^2\right){}_{k}}.\nonumber
 \end{align*}
 Thus the theorem is confirmed.
\qed
\end{pot222}
As mentioned earlier,  the  Bailey pair $(\alpha_n(y),\beta_n(y))$ given by \eqref{defnew-2} seems to have not appeared
in the present literature.
Based on  Theorem \ref{baileylemma}, it is easily found that
\begin{tl}\label{apartfrom} For $|q^2/ab|<1$, it holds
\begin{align}
\sum_{n=0}^{\infty}\frac{\poq{a,b}{n}}
{\poq{q^2/a,q^2/b}{n}}\left(\frac{1}{ab}\right)^nq^{n^2/2+5n/2}=\frac{\poq{q^2,q^2/ab}{\infty}}
{\poq{q^2/a,q^2/b}{\infty}}
{_3\phi_2}\bigg[\genfrac{}{}{0pt}{}{a,b,-q}{q^{3/2},-q^{3/2}};q,\frac{q^2}{ab}\bigg].\label{newnewnew}
 \end{align}
 In particular, we have
 \begin{align}
\sum_{n=0}^{\infty}
q^{3n^2/2+3n/2}=\poq{q}{\infty}\sum_{n=0}^{\infty}\frac{\poq{-q}{n}}{\poq{q}{n}\poqq{q}{n+1}}q^{n^2+n},\label{newnewnew-1}&\\
\sum_{n=0}^{\infty}\frac{\poq{a,-q^{3/2}}{n}}
{\poq{q^2/a,-q^{1/2}}{n}}\left(-\frac{1}{a}\right)^nq^{n^2/2+n}=\frac{\poq{q^2,q^{3/2}/a}{\infty}}{\poq{q^2/a,q^{3/2}}{\infty}}.\label{newnewnew-2}&
 \end{align}
\end{tl}
\pf Identity \eqref{newnewnew} is the special case $y=q$ of \eqref{slater-result}. Once letting  $a,b$ in
\eqref{newnewnew} tend to infinite, we obtain \eqref{newnewnew-1} while \eqref{newnewnew-2} results from setting
$b=-q^{3/2}$ and then applying the $q$-Gauss ${}_2\phi_1$ sum \cite[(II.8)]{10}.
\qed

It is interesting to note that  \eqref{newnewnew} can be used to an establishment of bilateral ${}_2\psi_2$ series
transformation. The reader may consult \cite[Eq.(5.1.1)]{10}  for the precise definition of the ${}_2\psi_2$ series.
\begin{tl} For $|ab|>1$, it holds
\begin{align}
\sum_{n=-\infty}^{\infty}\frac{\poq{aq,bq}{n}}
{\poq{q/a,q/b}{n}}\left(\frac{1}{ab}\right)^nq^{n^2/2+n/2}=2\frac{\poq{q^2,1/ab}{\infty}}
{\poq{q/a,q/b}{\infty}}{_3\phi_2}\bigg[\genfrac{}{}{0pt}{}{aq,bq,-q}{q^{3/2},-q^{3/2}};q,\frac{1}{ab}\bigg].\label{aaa}
 \end{align}
\end{tl}
\pf It suffices to calculate
\begin{align*}
\sum_{n=-\infty}^{\infty}\frac{\poq{aq,bq}{n}}
{\poq{q/a,q/b}{n}}\left(\frac{1}{ab}\right)^nq^{n^2/2+n/2}=\left\{\sum_{n=0}^{\infty}+\sum^{-1}_{n=-\infty}\right\}\frac{\poq{aq,bq}{n}}
{\poq{q/a,q/b}{n}}\left(\frac{1}{ab}\right)^nq^{n^2/2+n/2}
 \end{align*}
 while, in view of  \eqref{basic2},
\begin{align*}
\sum^{-1}_{n=-\infty}\frac{\poq{aq,bq}{n}}
{\poq{q/a,q/b}{n}}\left(\frac{1}{ab}\right)^nq^{n^2/2+n/2}&=\sum_{n=1}^{\infty}\frac{\poq{a,b}{n}}
{\poq{1/a,1/b}{n}}\left(\frac{1}{ab}\right)^nq^{n^2/2-n/2}\\
&\stackrel{n-1\to n}{=}\sum_{n=0}^{\infty}\frac{\poq{aq,bq}{n}}
{\poq{q/a,q/b}{n}}\left(\frac{1}{ab}\right)^nq^{n^2/2+n/2}.
 \end{align*}
 Therefore, we have
 \begin{align*}
\sum_{n=-\infty}^{\infty}\frac{\poq{aq,bq}{n}}
{\poq{q/a,q/b}{n}}\left(\frac{1}{ab}\right)^nq^{n^2/2+n/2}=2\sum_{n=0}^{\infty}\frac{\poq{aq,bq}{n}}
{\poq{q/a,q/b}{n}}\left(\frac{1}{ab}\right)^nq^{n^2/2+n/2}.
 \end{align*}
 Applying \eqref{newnewnew} under the replacements of $a,b$ with  $aq, bq$ to the above identity, we thereby obtain
 \eqref{aaa}.
\qed

Taking  \eqref{thetanew} into account, we easily find a short proof of Theorem \ref{mainthm-good}.
\begin{pot666}
It suffices to take $y=-a,-b$ in \eqref{thetanew} in succession, yielding
\begin{align*}
\theta(q,a)&=(q^2;q^2){}_\infty(a;q){}_{\infty}{}_{2}\phi_1\bigg[\genfrac{}{}{0pt}{}{-a/q,-a}{a^2};q^2,q^2\bigg],\\
\theta(q,b)&=(q^2;q^2){}_\infty(b;q){}_{\infty}{}_{2}\phi_1\bigg[\genfrac{}{}{0pt}{}{-b/q,-b}{b^2};q^2,q^2\bigg],
\end{align*}
where  the partial theta function $\theta(q,x)$ is defined to be the sum
 $$
\sum_{n=0}^\infty \tau(n)x^n.
$$
Multiplying these two identities together, we arrive at
\begin{align*}
\theta(q,a)\theta(q,b)
=(q^2;q^2){}_\infty^2(a,b;q){}_{\infty}{}_{2}\phi_1\bigg[\genfrac{}{}{0pt}{}{-a/q,-a}{a^2};q^2,q^2\bigg]
{}_{2}\phi_1\bigg[\genfrac{}{}{0pt}{}{-b/q,-b}{b^2};q^2,q^2\bigg].
\end{align*}
According to Theorem 1.1 of \cite{andrewswarnaar} by G.E. Andrews and S.O. Warnaar, we have
\begin{align}
\theta(q,a)\theta(q,b)=\poq{q,a,b}{\infty}\sum_{n=0}^\infty
\frac{\poq{ab/q}{2n}}{\poq{q,a,b,ab/q}{n}}q^n.\label{warnaar-andrews}
\end{align}
Both in together, after replacing $a,b$ with $-a,-b$ simultaneously and reformulating in terms of ${}_4\phi_3$ series,
leads us to the desired identity \eqref{thetanewidentity}. \qed
\end{pot666}

\subsection{Proof of Identity \eqref{4.1}}
A  careful look at \cite{andrews} shows that (\ref{4.1}) is given without direct proof. In what follows, we shall show
\eqref{4.1} in a different way from the above argument. This difference consists in that we utilize  a special Lagrange
inversion formula, instead of the transformation for ${}_3\phi_2$ series.
\begin{yl}[Lagrange inversion formula:\mbox{\cite[Example 2.2]{xrma}}] \label{fginversion} Let $\{x_n\}_{n\geq 0}$ be
any complex sequences such that $|1-x_nz|\neq 0$. If there exists the following expansion for certain analytic function
$F(z)$,
\begin{subequations}\label{4444}
 \begin{align}
 F(z)=\sum_{n=0}^\infty a_n\frac{z^n}{\prod_{i=1}^{n+1}(1-x_iz)},\label{4.11}
\end{align}
then the coefficients
\begin{align}
 a_n=[z^n] F(z)\prod_{i=1}^n(1-x_iz).\label{4.12}
\end{align}
\end{subequations}
\end{yl}
\begin{pot555}  By Lemma \ref{fginversion}, we first know that \eqref{4.1} is equivalent to
\begin{align}
\frac{q^{2n^2+2n+1}}{\poqq{q}{n+1}}=[z^n]\left(\sum_{k\geq
1}\frac{q^{k}z^{k-1}}{\poqq{q}{k}}\poqq{zq}{n}\right)\label{4.2},
\end{align}
or equivalently,
\begin{align}
\sum_{k=0}^{n}\bigg[\genfrac{}{}{0pt}{}{1+2n}{2k}\bigg]_q\poqq{q}{k}(-1)^kq^{k(k-1)}=q^{2n^2+n}.\label{known}
\end{align}
To make this more precise, we compute  the right-hand side of  \eqref{4.2} directly, namely
\begin{align*}
\mbox{RHS of
\eqref{4.2}}&=\sum_{k=0}^n\frac{q^{k+1}}{\poqq{q}{k+1}}[z^{n-k}]\frac{\poqq{zq}{\infty}}{\poqq{zq^{2n+1}}{\infty}}.
\end{align*}
Since
\begin{align*}
[z^{n-k}]\frac{\poq{zq^{1/2}}{\infty}}{\poq{zq^{n+1/2}}{\infty}}&=[z^{n-k}]
\sum_{i=0}^\infty\frac{\poq{q^{-n}}{i}}{\poq{q}{i}}(zq^{n+1/2})^i=\frac{\poq{q^{-n}}{n-k}}{\poq{q}{n-k}}q^{(n+1/2)(n-k)}\\
&=\frac{\poq{q^{-n}}{n}}{\poq{q}{n}}\frac{\poq{q^{-n}}{k}}{\poq{q}{k}}q^{(k+2
n^2+n)/2}=\bigg[\genfrac{}{}{0pt}{}{n}{k}\bigg]_q(-1)^{n-k}q^{(n-k)^2/2},
\end{align*}
wherein we have used the preceding relations \eqref{basic1} and \eqref{basic2}. On rescaling $q$ to $q^2$, it follows
that
\begin{align*}[z^{n-k}]\frac{\poqq{zq}{\infty}}{\poqq{zq^{2n+1}}{\infty}}&=\frac{\poqq{q^2}{n}}{\poqq{q^2}{k}\poqq{q^2}{n-k}}(-1)^{n-k}q^{(n-k)^2}.
\end{align*}
With these calculations, it is now easily verified that \eqref{4.2}
 is equivalent to
\begin{align*}
 q^{2n^2+2n+1}=\sum_{k=0}^nq^{k+1}\frac{\poqq{q}{n+1}\poqq{q^2}{n}}{\poqq{q}{k+1}\poqq{q^2}{k}}\frac{(-1)^{n-k}q^{(n-k)^2}}{\poqq{q^2}{n-k}}.
\end{align*}
Dividing by $q^{n+1}$ on both sides, we get
\begin{align*}
 q^{2n^2+n}&=\sum_{k=0}^nq^{k-n}\frac{\poq{q}{2n+1}}{\poq{q}{2k+1}\poq{q}{2(n-k)}}\frac{(-1)^{n-k}q^{(n-k)^2}\poq{q}{2(n-k)}}{\poqq{q^2}{n-k}}\\
&\stackrel{k\to
n-k}{=}\sum_{k=0}^n\bigg[\genfrac{}{}{0pt}{}{1+2n}{2k}\bigg]_q\frac{(-1)^{k}q^{k^2-k}\poq{q}{2k}}{\poqq{q^2}{k}}
=\sum_{k=0}^n\bigg[\genfrac{}{}{0pt}{}{1+2n}{2k}\bigg]_q\poqq{q}{k}(-1)^{k}q^{k^2-k}.
\end{align*}
Then  \eqref{known} follows. As such, all that we need to do is to check  the validity of   \eqref{known}, which is
asserted by   the $q$-Kummer (Bailey-Daum) sum \cite[(II.9)]{10}:
\begin{align}
\sum_{k=0}^{\infty}\frac{\poq{a,b}{k}}{\poq{q,aq/b}{k}}(-q/b)^k=\frac{\poq{-q}{\infty}\poqq{aq,aq^2/b^2}{\infty}}{\poq{-q/b,aq/b}{\infty}}.
\label{Bailey-Daum-old}
\end{align}
To be more precise, by setting $a=q^{-n}$ and letting  $b$ tend to infinity in \eqref{Bailey-Daum-old}, we first obtain
\begin{subequations}
\begin{align}
\sum_{k=0}^{n}\bigg[\genfrac{}{}{0pt}{}{n}{k}\bigg]_q\tau(n-k)Y(k)=X(n),
\label{Bailey-Daum}
\end{align}
where
\begin{align*}
 X(n)=\begin{cases}
  \displaystyle (-1)^kq^{k^2-k}\poqq{q}{k},\,\,\mbox{for}\,\, n=2k\\[2mm] \displaystyle
   0,\quad n=1
\pmod 2  \end{cases}
\quad\mbox{and}\qquad
 Y(n)=q^{n^2/2-n/2}.
\end{align*}
And then by inverting, it is easy to see that Identity \eqref{Bailey-Daum} is equivalent to
\begin{align}
\sum_{k=0}^{n}\bigg[\genfrac{}{}{0pt}{}{n}{k}\bigg]_qX(k)=Y(n),
\label{Bailey-Daum-inverting}
\end{align}
\end{subequations}
which, by inserting $X(n)$ and $Y(n)$ into, turns out to be  \eqref{known}. It gives the complete proof of
\eqref{4.1}.\qed
\end{pot555}
\section{Further discussions on two finite $q$-series summations}
In this section, we shall focus on two unusual finite sums closely related to the Andrews--Yee identities of Lemma
\ref{andrews-yee}.
Firstly, as Andrews and Yee demonstrated in \cite[Lemma 6]{andrews},  \eqref{generaltrans-7-77}
 is built on the key identity
\begin{equation}\sum_{k=0}^m\frac{\poq{q}{m+k}}{\poqq{q^2}{k}}q^k=\poqq{q^2}{m}.\label{generaltrans-888888-99}
\end{equation}
In the meantime, in the same paper \cite[Lemma 7]{andrews}  they recorded another finite $q$-series identity as below:
\begin{equation}\sum_{k=0}^m\frac{\poq{q}{m+k}}{\poqq{q^2}{k}}q^{2k}=\poqq{q}{m+1}+q^{m+1}\poqq{q^2}{m}.\label{generaltrans-888888-00}
\end{equation}
As for the importance of these two identities,  here we would like to quote  Andrew and Yee in \cite[p.3]{andrews} a comment on both
\eqref{generaltrans-888888-99} and \eqref{generaltrans-888888-00}:
\begin{quote}
 ``\emph{These results, while seemingly quite simple, are surprising for a couple of reasons. First, the sums do not
 terminate naturally, and second,we were unable to find these results in the $q$-series literature.}"
\end{quote}
It is this comment that  inspires us to investigate the sums
\begin{equation}
U_m(x):=\sum_{k=0}^m\bigg[\genfrac{}{}{0pt}{}{m+k}{k}\bigg]_q\frac{x^k}{\poq{-q}{k}}\label{generaltrans-888888-111}
\end{equation}
and
\begin{equation}
S_m(x,y):=\sum_{k=0}^m\bigg[\genfrac{}{}{0pt}{}{m+k}{k}\bigg]_q\frac{\poq{y}{k}}{\poq{x}{k}}q^k\label{generaltrans-888888-000}
\end{equation}
in a wider sense, in order to gain a better understanding how \eqref{generaltrans-888888-99} and
\eqref{generaltrans-888888-00} play a role in mock theta function identities.
\subsection{$q$-Difference equations for $U_m(x)$}
By considering $\{U_m(x)\}_{m\geq 0}$ as a polynomial sequence in $x$, we may set up two $q$-difference equations  for
it via the series rearrangement.
\begin{yl}\label{recurrenthem} Let $\{U_{m}(x)\}_{m\geq 0}$ be given by \eqref{generaltrans-888888-111}. Then for
integer $m\geq 0$, we have
\begin{align}
U_{m}(xq^2)&=U_{m}(x)-(1-q^{m+1})x U_{m+1}(x)+(x+1)x^{m+1}\frac{\poqq{q}{m+1}}{\poq{q}{m}},\label{ttt}\\
q^{m+1}U_{m}(xq)&=U_{m}(x)-(1-q^{m+1})U_{m+1}(x)+x^{m+1}\frac{\poqq{q}{m+1}}{\poq{q}{m}}.\label{ttt-1}
\end{align}
\end{yl}
\pf To show \eqref{ttt}, we start with \eqref{generaltrans-888888-111}. It is easy to verify that
\begin{align*}(1-q^{m})xU_{m}(x)&=\sum_{k=0}^m\frac{\poq{q^m}{k+1}}{\poqq{q^2}{k+1}}x^{k+1}(1-q^{2k+2})\nonumber\\
&=\sum_{k=1}^{m+1}\frac{\poq{q^m}{k}}{\poqq{q^2}{k}}x^{k}-\sum_{k=1}^{m+1}\frac{\poq{q^m}{k}}{\poqq{q^2}{k}}(xq^2)^{k}\\
&=U_{m-1}(x)-U_{m-1}(xq^2)+C_0(m,x),
\end{align*}
where, by some routine calculations,
\begin{align*}
C_0(m,x)=(x+1)x^{m}\frac{\poqq{q}{m}}{\poq{q}{m-1}}.
\end{align*}
Rearranging the last identity, we obtain
\begin{align*}U_{m-1}(xq^2)&=U_{m-1}(x)-(1-q^{m})xU_{m}(x)+(x+1)x^{m}\frac{\poqq{q}{m}}{\poq{q}{m-1}}.
\end{align*}
Shifting $m$ to $m+1$, we finally obtain the desired.

To establish \eqref{ttt-1}, let us write $S_m(x)=\poq{q}{m}U_m(x)$. Consider the difference
\begin{align}
S_m(x)-S_{m-1}(x)&=\sum_{k=0}^m\frac{\poq{q}{m+k}}{\poqq{q^2}{k}}x^k-\sum_{k=0}^{m-1}\frac{\poq{q}{m-1+k}}{\poqq{q^2}{k}}x^k\nonumber\\
&=\frac{\poq{q}{2m}}{\poqq{q^2}{m}}x^m+\sum_{k=0}^{m-1}\frac{\poq{q}{m+k}-\poq{q}{m-1+k}}{\poqq{q^2}{k}}x^k\nonumber\\
&=\poqq{q}{m}x^m-q^m\sum_{k=0}^{m-1}\frac{\poq{q}{m+k-1}}{\poqq{q^2}{k}}(xq)^k=\poqq{q}{m}x^m-q^mS_{m-1}(xq).
\label{ooo}
\end{align}
Replace $m$ with $m+1$ and recast \eqref{ooo} in terms of $U_m(x)$. We thus obtain the desired.
\qed
\begin{lz}
By virtue of \eqref{ttt-1} and the initial result \eqref{generaltrans-888888-99} (or see \cite[p.3]{andrews})
\begin{align*}
  U_m(q)=\poq{-q}{m}
\end{align*}
 we may find out all values of $U_m(q^{k}), k\geq 2.$ For instance, take $x=q$ in \eqref{ttt-1}. Then we have
\begin{align*}
q^{m+1}U_{m}(q^2)&=U_{m}(q)-(1-q^{m+1})U_{m+1}(q)+q^{m+1}\frac{\poqq{q}{m+1}}{\poq{q}{m}}\\
&=\poq{-q}{m}-(1-q^{m+1}) \poq{-q}{m+1}+q^{m+1}\frac{\poqq{q}{m+1}}{\poq{q}{m}}\\
&=\poq{-q}{m}q^{2m+2}+q^{m+1}\frac{\poqq{q}{m+1}}{\poq{q}{m}}.
\end{align*}
It amounts to \eqref{generaltrans-888888-00}.
\end{lz}
With Lemma \ref{recurrenthem} in hand, it is easily found that $U_m(x)$ satisfies a second order  $q$-difference
equation.
\begin{dl} Let $\{U_{m}(x)\}_{m\geq 0}$ be given by \eqref{generaltrans-888888-111}. Then for integer $m\geq 0$, it
holds
\begin{align}U_{m}(xq^2)=(1-x)U_m(x)+xq^{m+1}U_m(xq)+x^{m+1}\frac{\poqq{q}{m+1}}{\poq{q}{m}}.
\end{align}
\end{dl}
\pf It only needs to restate \eqref{ttt-1} in the form
\begin{align}
(1-q^{m+1})U_{m+1}(x)=U_{m}(x)-q^{m+1}U_{m}(xq)+x^{m+1}\frac{\poqq{q}{m+1}}{\poq{q}{m}}.\label{add}
\end{align}
Substituting \eqref{add} into \eqref{ttt}, then we have
\begin{align*}
U_{m}(xq^2)&=U_{m}(x)-x\left(U_{m}(x)-q^{m+1}U_{m}(xq)+x^{m+1}\frac{\poqq{q}{m+1}}{\poq{q}{m}}\right)+(x+1)x^{m+1}\frac{\poqq{q}{m+1}}{\poq{q}{m}}\\
&=(1-x)U_{m}(x)+xq^{m+1}U_{m}(xq)+x^{m+1}\frac{\poqq{q}{m+1}}{\poq{q}{m}}.
\end{align*}
\qed

As another application of Lemma \ref{recurrenthem}, we can deduce a three-term recurrence relation for
$\{U_m(x)\}_{m\geq 0}$ without the $q$-difference operator involved.
\begin{dl} Let $\{U_{m}(x)\}_{m\geq 0}$ be given by \eqref{generaltrans-888888-111}. Then for integer $m\geq 0$, it
holds
\begin{align}(1-q^{m+2})U_{m+2}(x)+(x q^{2
m+3}-q-1)U_{m+1}(x)+q(1+q^{m+1})U_{m}(x)=(x-q)\frac{\poq{q^{m+2}}{m}}{\poqq{q^2}{m}}x^{m+1}.\label{qqq}
\end{align}
\end{dl}
\pf It suffices to compute, by using  \eqref{ttt-1},  that
\begin{align}
&q^{2m+3}U_{m}(xq^2)=q^{m+2}U_{m}(xq)-(1-q^{m+1})q^{m+2}U_{m+1}(xq)+x^{m+1}q^{2m+3}\frac{\poqq{q}{m+1}}{\poq{q}{m}}.\nonumber
\end{align}
Next, applying  \eqref{ttt-1} again, we obtain
\begin{align}
q^{2m+3}U_{m}(xq^2)&=qU_m(x)-(1+q)(1-q^{m+1})U_{m+1}(x)+(1-q^{m+1})(1-q^{m+2})U_{m+2}(x)\nonumber\\
&+qx^{m+1}\frac{\poqq{q}{m+1}}{\poq{q}{m}}-x^{m+2}\frac{\poqq{q}{m+2}}{\poq{q}{m}}+x^{m+1}q^{2m+3}\frac{\poqq{q}{m+1}}{\poq{q}{m}}.\label{leisire}
\end{align}
Lastly, multiplying  \eqref{ttt} by $q^{2m+3}$ and subtracting from \eqref{leisire}, we conclude that
\begin{align*}
0&=q(1-q^{2m+2})U_m(x)+(xq^{2m+3}-1-q)(1-q^{m+1})U_{m+1}(x)\\
&+(1-q^{m+1})(1-q^{m+2})U_{m+2}(x)+qx^{m+1}\frac{\poqq{q}{m+1}}{\poq{q}{m}}-x^{m+2}\frac{\poqq{q}{m+1}}{\poq{q}{m}}.
\end{align*}
This gives the expression of $U_{m+2}(x)$ in terms of $U_{m}(x)$ and $U_{m+1}(x)$.
\qed
\begin{remark}
Interestingly, in light of the WZ method, we also know that the summand of $U_m(x)$, denoted by $V_m(k;x)$,   satisfies
the recurrence relation
\begin{align*}V_{m+2}(k;x)&+\frac{-x q^{2
m+3}+q+1}{q^{m+2}-1}V_{m+1}(k;x)+\frac{q(q^{m+1}+1)}{1-q^{m+2}}V_{m}(k;x)\nonumber\\
&=V_m(k;x)R(m,k;x)-V_m(k-1;x)R(m,k-1;x),
\end{align*}
where the certification
$$
R(m,k;x)=\frac{x q^{2 m+3} (1-q^{m+k+1})}{(1-q^{m+1})(1-q^{m+2})}.$$
\end{remark}

\subsection{A general transformation for $S_m(x,y)$ and its implications}
In this part, we turn to  possible transformations associated with the sequence $\{S_m(x,y)\}_{m\geq 0}$.
\begin{dl}\label{generaltrans} For any $q:|q|<1$ and integer $m\geq 0$, it always holds
\begin{align}\displaystyle\frac{
\poq{yq^2/x}{m}}{\poq{q^2/x}{m}}\sum_{k=0}^m\bigg[\genfrac{}{}{0pt}{}{m+k}{k}\bigg]_q\frac{\poq{y}{k}}{\poq{x}{k}}q^k
=1&+\frac{q+x}{x-yq}\sum _{k=1}^m
\genfrac{[}{]}{0pt}{}{2k-1}{k}_{q}\frac{\poq{y,yq/x}{k}}{\poq{x,q^2/x}{k}}q^k\nonumber\\
&-\frac{yq}{x-yq}\sum _{k=1}^m
\genfrac{[}{]}{0pt}{}{2k}{k}_{q}\frac{\poq{y,yq/x}{k}}{\poq{x,q^2/x}{k}}q^{2k}.\label{masterid}
\end{align}\end{dl}
\pf Performing as before, we first show by  the WZ method  that  $\{S_m(x,y)\}_{m\geq 0}$ satisfies the first order
recurrence relation
\begin{align}S_{m}(x,y)&-\frac{x-q^{m+1}}{x-yq^{m+1}}S_{m-1}(x,y)
\nonumber\\
&=q^{m}\frac{q+x-q^{m+1}(1+q^m)y}{x-q^{m+1}y}\bigg[\genfrac{}{}{0pt}{}{2m-1}{m}\bigg]_q\frac{\poq{y}{m}}{\poq{x}{m}}.\label{ppp}
\end{align}
Next,  we solve this  recurrence relation for $S_{m}(x,y)$ and obtain
\begin{align*}
S_m(x,y)=\frac{\poq{q^2/x}{m}}{ \poq{yq^2/x}{m}}
\bigg(1+\frac{1}{x-yq}\sum _{k=1}^m q^k
(q+x-y(q^k+1)q^{k+1})\genfrac{[}{]}{0pt}{}{2k-1}{k}_{q}\frac{\poq{y,yq/x}{k}}{\poq{x,q^2/x}{k}}\bigg).
\end{align*}
A combination of this expression and the definition  \eqref{generaltrans-888888-000}, after a bit simplification, leads
us to \eqref{masterid}.
\qed

It is worth mentioning, as applications of Theorem \ref{generaltrans},  that \eqref{masterid} implies several noteworthy
transformations.
\begin{tl}For arbitrary integer $m\geq 0$, we have
 \begin{align}\frac{\poq{-yq}{m}}{\poq{-q}{m}}
\sum_{k=0}^m\frac{\poq{q^{m+1},y}{k}}{\poqq{q^2}{k}}q^k
&=1+\frac{y}{1+y}
\sum_{k=1}^m \genfrac{[}{]}{0pt}{}{2k}{k}_q
\frac{\poqq{y^2}{k}}{\poq{-q}{k}^2}q^{2k},\label{special-2}\\
&\nonumber\\
\frac{1}{\poq{q^2/x}{m}}\sum_{k=0}^m\bigg[\genfrac{}{}{0pt}{}{m+k}{k}\bigg]_q\frac{q^k}{\poq{x}{k}}&=1+(1+q/x)\sum
_{k=1}^m \genfrac{[}{]}{0pt}{}{2k-1}{k}_{q}\frac{q^k}{\poq{x,q^2/x}{k}}.\label{zero}
\end{align}
\end{tl}
\pf
 In fact,  \eqref{special-2} is the consequence of \eqref{masterid} for $x=-q$ while \eqref{zero} is the special case
 $y=0$ of \eqref{masterid}.
\qed

\begin{tl}For arbitrary integer $m\geq 0$, we have
\begin{align}\poq{x}{m+1}\sum_{k=0}^m\bigg[\genfrac{}{}{0pt}{}{m+k}{k}\bigg]_q\frac{q^k}{\poq{xq}{k}}
+\big(1/x;q\big)_{m+1}\sum_{k=0}^m\bigg[\genfrac{}{}{0pt}{}{m+k}{k}\bigg]_q\frac{q^k}{\poq{q/x}{k}}=\big(x,1/x;q\big)_{m+1}.\label{3.6}
\end{align}
\end{tl}
\pf It suffices to replace  $x$ with $xq$ in \eqref{zero}. We have
\begin{align}\frac{1}{\poq{q/x}{m}}\sum_{k=0}^m\bigg[\genfrac{}{}{0pt}{}{m+k}{k}\bigg]_q\frac{q^k}{\poq{xq}{k}}=1+(1+1/x)\sum
_{k=1}^m \genfrac{[}{]}{0pt}{}{2k-1}{k}_{q}\frac{q^k}{\poq{xq,q/x}{k}}.\label{one}
\end{align}
Next, changing  $x$ to $1/x$ in \eqref{one}, we  arrive at
\begin{align}\frac{1}{\poq{xq}{m}}\sum_{k=0}^m\bigg[\genfrac{}{}{0pt}{}{m+k}{k}\bigg]_q\frac{q^k}{\poq{q/x}{k}}=1+(1+x)\sum
_{k=1}^m \genfrac{[}{]}{0pt}{}{2k-1}{k}_{q}\frac{q^k}{\poq{xq,q/x}{k}}.\label{two}
\end{align}
Multiplying $\eqref{one}$ by $x$ and then subtracting \eqref{two} from  the resulting identity,  we get
\begin{align}\frac{x}{\poq{q/x}{m}}\sum_{k=0}^m\bigg[\genfrac{}{}{0pt}{}{m+k}{k}\bigg]_q\frac{q^k}{\poq{xq}{k}}
-\frac{1}{\poq{xq}{m}}\sum_{k=0}^m\bigg[\genfrac{}{}{0pt}{}{m+k}{k}\bigg]_q\frac{q^k}{\poq{q/x}{k}}=x-1.\label{three}
\end{align}
Lastly, dividing by $x-1$ on both sides of \eqref{three} and multiplying the resulted with $\poq{x,1/x}{m+1}$ gives the
desired identity.
\qed

It is worth pointing out that from Identity \eqref{3.6} we can recover  the classical $q$-Chu-Vandermonde identity
\cite[(II.7)]{10} and the well-known Jacobi triple product identity \cite[(II.28)]{10}.
\begin{lz}
\begin{description}
    \item [(i)]\,\,\mbox{($q$-Chu-Vandermonde identity)} For any integer $r\geq 0$, it holds
\begin{align}
\sum_{k=0}^m\bigg[\genfrac{}{}{0pt}{}{m+k}{m}\bigg]_q\bigg[\genfrac{}{}{0pt}{}{m-k}{r}\bigg]_qq^{(r+1)k}
=\bigg[\genfrac{}{}{0pt}{}{2m+1}{m+r+1}\bigg]_q.\label{chu-van}
\end{align}
  \item [(ii)]\,\mbox{(Jacobi triple product identity)}
  \begin{align}
\big(x,q/x,q;q\big)_{\infty}=\sum_{n=-\infty}^{\infty}(-1)^nq^{n(n-1)/2}x^n.\label{triple}
\end{align}
\end{description}
\end{lz}
\pf (i) Observe that
\[\poq{x,1/x}{m+1}=(1-1/x)\tau(m)(q/x)^m\poq{q^{-m}x}{2m+1}.\]
After dividing  both sides of \eqref{3.6} by $1-1/x$, we obtain
\begin{align*}-x\poq{xq}{m}\sum_{k=0}^m\bigg[\genfrac{}{}{0pt}{}{m+k}{k}\bigg]_q\frac{q^k}{\poq{xq}{k}}
+\big(q/x;q\big)_{m}\sum_{k=0}^m\bigg[\genfrac{}{}{0pt}{}{m+k}{k}\bigg]_q\frac{q^k}{\poq{q/x}{k}}=\tau(m)(q/x)^m\poq{q^{-m}x}{2m+1},
\end{align*}
or equivalently,
\begin{align*}-x\sum_{k=0}^m\bigg[\genfrac{}{}{0pt}{}{m+k}{k}\bigg]_qq^k\poq{xq^{k+1}}{m-k}
+\sum_{k=0}^m\bigg[\genfrac{}{}{0pt}{}{m+k}{k}\bigg]_qq^k\poq{q^{k+1}/x}{m-k}=\tau(m)(q/x)^m\poq{q^{-m}x}{2m+1}.
\end{align*}
By equating the coefficients of $x^{r+1}$, the claimed follows.

(ii)  To establish \eqref{triple}, we first let $m$ tend to infinity in \eqref{3.6}.  Then it follows
\begin{align*}\poq{x}{\infty}\sum_{k=0}^{\infty}\frac{q^k}{\poq{q,xq}{k}}
+\big(1/x;q\big)_{\infty}\sum_{k=0}^{\infty}\frac{q^k}{\poq{q,q/x}{k}}=\big(x,1/x;q\big)_{\infty}.
\end{align*}
Multiplying both sides with $\poq{q}{\infty}/(1-1/x)$,  we therefore obtain
\begin{align} -x~\poq{q,xq}{\infty}\sum_{k=0}^{\infty}\frac{q^k}{\poq{q,xq}{k}}
+\big(q,q/x;q\big)_{\infty}\sum_{k=0}^{\infty}\frac{q^k}{\poq{q,q/x}{k}}=\big(x,q/x,q;q\big)_{\infty}.\label{none}
\end{align}
Observe that  \eqref{warnaar-andrews} (also see \cite[(2.1a)]{andrewswarnaar}) implies
\begin{align*}
\sum_{k=0}^{\infty} \tau(k) x^k=\poq{q,x}{\infty}\sum_{k=0}^{\infty} \frac{q^k}{\poq{q,x}{k}}.
\end{align*}
This allows us to reformulate  \eqref{none} in the form
\begin{align*}
-x\sum_{k=0}^{\infty} \tau(k) (xq)^k
+\sum_{k=0}^{\infty} \tau(k) (q/x)^k=\big(x,q/x,q;q\big)_{\infty}.
\end{align*}
After a bit of manipulation, \eqref{triple} follows.
\qed

\begin{tl} For arbitrary integer $m\geq 0$,
\begin{align}
\poq{x}{m+1}\sum_{k=0}^m\bigg[\genfrac{}{}{0pt}{}{m+k}{k}\bigg]_q\frac{q^k}{\poq{xq}{k}}
=\bigg[\genfrac{}{}{0pt}{}{2m+1}{m}\bigg]_q+\frac{1}{1+q^{m+1}}\sum_{k=1}^{m+1}\tau(k)(1+q^k){2m+2\brack
m+1-k}_qx^k.\label{3.600}
\end{align}
\end{tl}
\pf Since,  easily proved by induction, that
\begin{eqnarray*}
(x,1/x;q)_n=\frac{1}{1+q^n}\sum_{k=-n}^n\tau(k)(1+q^k){2n\brack n-k}_qx^k,
\end{eqnarray*}
thus \eqref{3.6} is equivalent to
\begin{align*}\poq{x}{m+1}\sum_{k=0}^m\bigg[\genfrac{}{}{0pt}{}{m+k}{k}\bigg]_q\frac{q^k}{\poq{xq}{k}}
&+\big(1/x;q\big)_{m+1}\sum_{k=0}^m\bigg[\genfrac{}{}{0pt}{}{m+k}{k}\bigg]_q\frac{q^k}{\poq{q/x}{k}}\nonumber\\
&=\frac{1}{1+q^{m+1}}\sum_{k=-m-1}^{m+1}\tau(k)(1+q^k){2m+2\brack m+1-k}_qx^k.
\end{align*}
Then we separate the sums into two sums on both sides according as the index $k$ of $x^k$ is nonnegative or not. By
comparing coefficients, we have
\begin{align}\poq{x}{m+1}\sum_{k=0}^m\bigg[\genfrac{}{}{0pt}{}{m+k}{k}\bigg]_q\frac{q^k}{\poq{xq}{k}}+\sum_{k=0}^m\bigg[\genfrac{}{}{0pt}{}{m+k}{k}\bigg]_qq^k
=\frac{1}{1+q^{m+1}}\sum_{k=0}^{m+1}\tau(k)(1+q^k){2m+2\brack m+1-k}_qx^k.\label{3.6001}
\end{align}
Note that  the special case $r=0$ of \eqref{chu-van} yields
\begin{align*}\sum_{k=0}^m\bigg[\genfrac{}{}{0pt}{}{m+k}{k}\bigg]_qq^k=\bigg[\genfrac{}{}{0pt}{}{2m+1}{m}\bigg]_q.\label{qqq1}
\end{align*}
This simplifies \eqref{3.6001} to the desired identity.
\qed

\section{Concluding remarks}
We conclude our paper with a few remarks on further study of $\omega(z;q)$ and $v(z;q)$.  For this, let us consider the
limiting case of \eqref{generaltrans-new} as n tends to infinite and  $c=q, e=zq$, stated as follows.
\begin{tl} For $|dz/b|,|z|<1$, it holds
\begin{eqnarray}
\sum_{k=0}^{\infty}\tau(k)
\frac{\poq{b}{k}}{\poq{d,zq}{k}}\bigg(\frac{dz}{b}\bigg)^k=
(1-z)\sum_{k=0}^{\infty}\frac{\poq{d/b}{k}}{\poq{d}{k}}z^k.
\label{sect4-generaltrans}
\end{eqnarray}
\end{tl}
From this transformation it is easy to deduce
\begin{dl}\label{th41} For $|y|<1,|z|<1$, we have
\begin{eqnarray}
\sum_{k=0}^{\infty}
\frac{q^{k^2}(yz)^k}{\poq{y,z}{k+1}}&=&\sum_{k=0}^\infty\frac{z^k}{\poq{y}{k+1}}=
\sum_{k=0}^\infty\frac{y^k}{\poq{z}{k+1}},\label{generaltrans-4}\\
\sum_{k=0}^\infty\frac{\tau(k)(yz)^k}{\poq{z}{k+1}}&=&\sum_{k=0}^\infty\poq{y}{k}z^k.\label{generaltrans-6}
\end{eqnarray}
\end{dl}
\pf (i) Identity \eqref{generaltrans-4} follows immediately by taking $d=yq$ and letting $b\to\infty$ in
\eqref{sect4-generaltrans}.
(ii)  \eqref{generaltrans-6} is a direct result while we set $d/b=y$ and then put $b=0$ in \eqref{sect4-generaltrans}.
\qed

\begin{remark}Define
\begin{align}
\omega_q(y,z):=
\sum_{k=0}^{\infty}
\frac{q^{k^2}(yz)^k}{\poq{y,z}{k+1}},\quad v_q(y,z):=
\sum_{k=0}^{\infty}\frac{\tau(k)(yz)^k}{\poq{z}{k+1}}.
\end{align}
Then it becomes very clear that
\begin{align}
\omega(z;q)=\omega_{q^2}(q,zq),\quad v(z;q)=v_{q^2}(q/z,-zq).
\end{align}
where $\omega(z;q)$ and $v(z;q)$ are the same as before. In particular, \eqref{generaltrans-6} under the parametric
replacement $(q,y,z)\to (q^2,zq,-q)$ is just  Theorem 2 of \cite{andrews}. Some recent results on $\omega_q(y,z)$ have
been made by Y.S. Choi in \cite[Section 4]{newadded}.
\end{remark}
Based on such supporting evidence,  an interesting problem arises from Theorem \ref{th41} as following.
\begin{wenti} Given $\omega_q(y,z)$ and $v_q(y,z)$, are we able to find certain results   analogue to  Theorems 1 and 2
of \cite{andrews} for $\omega(z;q)$ and $v(z;q)$?
\end{wenti}

\bibliographystyle{amsplain}

\end{document}